\documentclass[twocolumn]{revtex4}
\usepackage{graphicx}

\begin{document}

\title{Variational optimization of probability measure spaces resolves \\ the chain store paradox}

\author{Michael J. Gagen}

\affiliation{Institute for Molecular Bioscience, University of
Queensland, Brisbane, Qld 4072, Australia}

\email{m.gagen@imb.uq.edu.au}

\author{Kae Nemoto}

\affiliation{National Institute of Informatics, 2-1-2
Hitotsubashi, Chiyoda-ku, Tokyo 101-0843, Japan}

\email{nemoto@nii.ac.jp}

\date{26 April 2006}

\begin{abstract}
In game theory, players have continuous expected payoff functions
and can use fixed point theorems to locate equilibria. This
optimization method requires that players  adopt a particular type
of probability measure space. Here, we introduce alternate
probability measure spaces altering the dimensionality,
continuity, and differentiability properties of what are now the
game's expected payoff functionals. Optimizing such functionals
requires generalized variational and functional optimization
methods to locate novel equilibria.  These variational methods can
reconcile game theoretic prediction and observed human behaviours,
as we illustrate by resolving the chain store paradox. Our
generalized optimization analysis has significant implications for
economics, artificial intelligence, complex system theory,
neurobiology, and biological evolution and development.
\end{abstract}

\maketitle

\section{Introduction}

In game theory, as formalized by von Neumann and Morgenstern
\cite{vonNeumann_44}, Nash \cite{Nash_50_48,Nash_51_28}, and Kuhn
\cite{Kuhn_1953}, rational players with common knowledge of
rationality (CKR) locate equilibria by using fixed point theorems
to optimize continuous expected payoff functions. These expected
payoff functions, according to probability measure theory
\cite{Bauer_1981,Pfeiffer_1990,Kelly_94}, can only be defined
after the adoption of a suitable probability measure space
supporting appropriate random variables, functions, and
probability distributions. For instance, mixed strategy
probability measure spaces were used by von Neumann and
Morgenstern \cite{vonNeumann_44} and Nash
\cite{Nash_50_48,Nash_51_28}, while behavioural strategy
probability measure spaces were introduced by Kuhn
\cite{Kuhn_1953}.  In addition, correlated strategy probability
measure spaces were introduced by Aumann to model communication
channels between players \cite{Aumann_74_67}. In this last case,
communications necessitate a change of probability measure space,
however a change of probability space does not always require
communication.  Consequently, in this paper we introduce a method
to analyze games using the infinite number of different
probability measure spaces available to describe any given game
and set of expected payoffs
\cite{Bauer_1981,Pfeiffer_1990,Kelly_94}. Our particular interest
lies in the class of probability measure spaces which is
consistent with the given game information constraints. That is,
we consider only probability measure spaces which are consistent
with rationality, CKR, and no communication channels between
players. Such probability measure spaces can exist, as we show
later, simply because a number of different probability measure
spaces are consistent with information flow via the game history
set without any communication channels. In this paper, we suppose
players may freely alter their choice of probability measure space
among all those consistent with no communications or any other
alteration in the game, in contrast to, for instance, previous
work on correlated equilibria \cite{Aumann_74_67}. For many games,
a change in the underlying probability measure space will not
affect equilibria---witness the equivalence of mixed and
behavioural strategies in games of perfect recall
\cite{Kuhn_1953}. However, in this paper, we argue that there
exist games in which altering the choice of probability measure
space will alter strategic equilibria. Assuming rationality, CKR,
and the usual game information constraints, players can search an
enlarged space of alternate probability measure spaces to optimize
their expected payoffs, and thereby locate novel equilibria
improving their outcomes over those achieved using only the
conventional mixed or behavioural strategy probability spaces of
game theory.

In this paper, we assess for the first time whether the set of
equilibria of any arbitrary game are entirely invariant under the
altered mathematical parameterizations defined by different
probability measure spaces. It does appear that equilibria are
indeed invariant under alternate probability measure spaces for
single-player and multiple-player-single-stage games. However,
equilibria are not invariant under altered choice of probability
measure space for multiple-player-multiple-stage games.  In these
games, the adoption of alternate probability measure spaces by
players can so alter the parameterized expected payoff functions
as to generate entirely novel sets of equilibria.

Demonstrating this requires a significant generalization of the
usual optimization methods of game theory.  This is because
alternate probability measure spaces and parameterizations can
alter the functional form, dimensionality, continuity and
differentiability properties of what must now be treated as
expected payoff functionals (not functions).  As a result, the
multiple-player calculus methods (essentially fixed point
theorems) suitable for expected payoff functions defined over
continuous probability simplexes are insufficient. To optimize
expected payoff functionals, we must generalize the variational
and functional optimization techniques used in, for instance,
general equilibrium and Cass-Koopmans style optimal growth
analysis \cite{Cass_1965_233,Koopmans_1965_1}, Ramsey-style
multiple stage optimization
\cite{Ramsey_1928_543,Kamien_1991,Chiang_2000}, and continuous
time differential games \cite{Dockner_2000}.  Suitably
generalized, these variational and functional optimization
techniques can reconcile game theoretic prediction and observed
human behaviour as we illustrate using Selton's chain store
paradox \cite{Selten_78_12}.  In this game, backwards induction
predicts that a monopolist never fights new market entrants even
though, in practice, most monopolists will indeed fight new
entrants and thereby improve their payoffs. This led Selton to
conclude ``mathematically trained persons recognize the logical
validity of the induction argument, but they refuse to accept it
as a guide to practical behavior." \cite{Selten_78_12}.  This
stark contrast makes this game a suitable vehicle for the
presentation of our new methods.

\section{Variational optimization of probability measure spaces}

We consider the general strategic optimization problem faced by
two players $X$ and $Y$ seeking to maximize their respected
expected payoffs $\langle\Pi^X\rangle$ and $\langle\Pi^Y\rangle$
in a game where $X$ chooses events $x$ and $Y$ chooses events $y$
to generate respective payoff outcomes for each player of
$\Pi^{X}(x,y)$ and $\Pi^{Y}(x,y)$.  The chosen events $x\times y$
are contained in $\Omega^X\times\Omega^Y$, the set of all possible
events in the game and in both player's chosen ``roulette"
randomization devices. These devices are used by players to avoid
their choices being forecast and exploited, with the result that
the choice of events is described using a joint probability
distribution $P^{XY}_{xy}$. As is required in probability measure
theory \cite{Bauer_1981,Pfeiffer_1990,Kelly_94}, the definition of
this joint probability distribution requires player $X$ to adopt a
probability measure space ${\cal P}^X$, and player $Y$ to adopt a
probability measure space ${\cal P}^Y$, such that the joint
product probability measure space ${\cal P}^X\times{\cal P}^Y$
supports the probability measure $P^{XY}_{xy}$. We allow players
to vary their choice of probability measure space to maximize
their expected payoffs. Altogether, the strategic optimization
problem facing each player is
\begin{eqnarray}
 X\!\!: \; \max_{{\cal P}^X}\;\; \langle \Pi^X \rangle \!\!
    &=& \!\! \int_{\Omega^X\!\!\times\Omega^Y}
         dP^{XY}_{xy} \;\;
         \Pi^X(x,y) \nonumber \\
    & & \\
 Y\!\!: \; \max_{{\cal P}^Y} \;\; \langle \Pi^Y \rangle \!\!
    &=& \!\! \int_{\Omega^X\!\!\times\Omega^Y}
         dP^{XY}_{xy} \;\;
         \Pi^Y(x,y). \nonumber
\end{eqnarray}
Here, expected payoffs for each player $Z\in\{X,Y\}$ are defined
by a Lebesgue integral over all possible game and roulette events
$\Omega^X\times\Omega^Y$ of payoffs $\Pi^Z(x,y)$ resulting from
particular game events $(x,y)$ weighted by the joint probability
measure of those events occurring $P^{XY}_{xy}$. The optimization
involves each player $Z$ maximizing their expected payoff over
every possible joint probability measure space that might be
adopted ${\cal P}^X\times{\cal P}^Y$, where ${\cal
P}^Z=\{\Omega^Z,\sigma^Z,P^Z\}$ is defined in terms of an
appropriate event set $\Omega^Z$ modelling all game and roulette
device events, a suitable sigma-algebra $\sigma^Z$, and an
appropriate probability measure $P^Z$.

Game theory has not previously allowed rational players to vary
their choice of probability space to maximize their expected
payoffs. This is largely because von Neumann and Morgenstern's
original goal was to formulate strategic plans assessing every
possible move in a game \cite{vonNeumann_44}, and they considered
this goal required only that each player adopt a particular
probability measure space defining mixed strategies in any game.
(Kuhn later introduced alternate behavioural strategy probability
measure spaces providing an equivalent analysis in games of
perfect recall \cite{Kuhn_1953}.) While never stated explicitly,
this restriction essentially limits the search space of the
players so they can only optimize over the probability parameters
of a single type of probability space using fixed point theorems
to locate Nash equilibria. In contrast, we argue that, under CKR,
players can search every alternate probability space consistent
with game information constraints by using generalized variational
and functional optimization techniques. In the remainder of this
section, we seek to explain heuristically why such a generalized
analysis can generate novel and improved equilibria, and thus
reconcile game theoretic prediction and observed human behaviours.

Alternate probability measure spaces can support different
equilibria in strategic situations as each adopted probability
space can mathematically parameterize the same random event in
very different ways. For example, consider a player $X$ seeking to
optimize a binary outcome specified by a random variable taking
value $x=0$ with probability $P^X(0)$ or $x=1$ with probability
$P^X(1)$. These probabilities can be characterized in terms of a
single probability parameter $p$ by tossing a biased coin, or in
terms of five probability parameters $(p_1,p_2,p_3.p_4,p_5)$ say
by using a biased dice. An alternative probability measure space
might employ two sequentially tossed, independent, biased coins
producing outcomes $u=1$ with probability $p$, while if $u=0$ then
$v=1$ with probability $q$ and if $u=1$ then $v=1$ with
probability $r$. The subsequent adoption of the random variable
$x=\delta_{u1}\delta_{v1}$ defines $P^X(1)=P(u=1,v=1)=pr$. (Here,
$\delta_{ab}=1$ if $a=b$ and zero otherwise.) As a last
illustration, consider a probability measure space in which the
above two biased coins are now perfectly correlated via
$P(u,v)=P(u)P(v|u)=P(u)\delta_{uv}=P(u)$. In this case, the known
perfect correlation introduces a delta function to reduce the
dimensionality of the joint distribution $P(u,v)$ giving
$P^X(1)=P(u=1,v=1)=p$.  In general, when parameterized using
different probability measure spaces, a given probability
possesses alternate functional forms with different
dimensionality, correlation, continuity, and differentiability
properties.

This changeability of functional form and dimensionality requires
generalized variational and functional optimization methods be
used to optimize strategic decisions.  The generalized methods we
develop extend the calculus of variations which typically
optimizes a functional $F[f(x),f'(x)]$ of known form, and where
the functional $F$, the function $f(x)$, and the gradient $f'(x)$
have specified differentiability properties. For instance, a
shortest path problem seeks to optimize the known functional
$F[f(x),f'(x)]=\sqrt{1+f'^2}$ via
\begin{equation}
    \max_{f,f'} \;\; I = \int_a^b \sqrt{1+f'^2} dx.
\end{equation}
Similarly, the shortest time or Brachistochrone problem optimizes
the known functional $F[f(x),f'(x)]=\sqrt{\frac{1+f'^2}{2gf}}$ via
\begin{equation}
    \max_{f,f'} \;\; I = \int_a^b \sqrt{\frac{1+f'^2}{2gf}} dx.
\end{equation}
Lastly, a typical multiple stage Ramsey-style utility maximization
problem optimizes
\begin{equation}
    \max_{f,f'} \;\; I = \int_a^b e^{-rx} F[f-f'] dx,
\end{equation}
where now only the functional dependencies and certain
differentiability properties of the functional $F[f(x),f'(x)]$ are
specified. To our knowledge, all applications of the calculus of
variations place severe restrictions on the range of variation of
the form of the functional being optimized, so much so that a
problem with an entirely arbitrary functional would be considered
ill defined. In contrast, in a strategic optimization problem,
players are able to arbitrarily vary their choice of probability
measure space to alter all of the functional form, the
dimensionality, and the continuity and differentiability
properties of the functional being optimized.  Heuristically, in
single player terms, the optimization problem becomes
\begin{equation}
    \max_{f,f'} \;\; I =
      \left\{
 \begin{array}{c}
   \vdots \\
    \\
   \int_a^b \sqrt{1+f'^2} dx \\
    \\
   \int_a^b \sqrt{\frac{1+f'^2}{2gf}} dx \\
    \\
   \int_a^b e^{-rx} F[f-f'] dx \\
    \\
   \vdots \\
 \end{array}
      \right. .
\end{equation}
That is, each player has the option of first choosing a
parameterizing probability measure space to alter the functional
form, dimensionality, continuity and differentiability properties
of the functionals being optimized, and only then to optimize the
chosen functional over all possible variations of $f(x)$ and
$f'(x)$. More importantly, each of their choices affects their
opponent's functionals, while at the same time, their opponent's
decisions are similarly altering their own functionals.

\begin{figure}[htb]
\centering
\includegraphics[width=\columnwidth,clip]{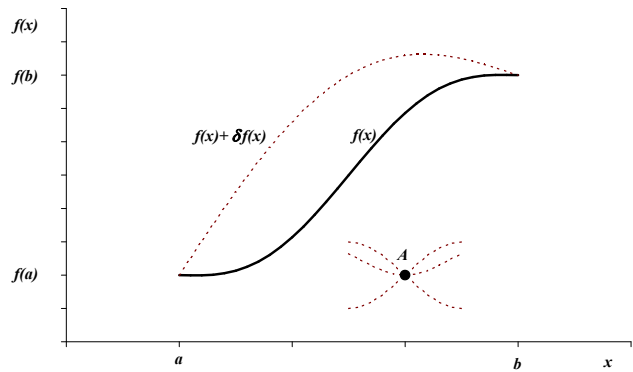}
\caption{\em  The variational optimization of the functional
$F[f(x),f'(x)]$ requires the variation of both the function
$f(x)\rightarrow f(x)+\delta f(x)$ and independently, its gradient
$f'(x)\rightarrow f'(x)+\delta f'(x)$ over the region $[a,b]$.
That is, through any point ``$A$", every possible gradient must be
considered in order to perform a complete variational analysis.
 \label{f_variational_calculus}}
\end{figure}

We suggest that this variability of the strategic functionals
means that optimization requires independent examination of every
possible functional, and every possible functional gradient, that
might be defined by the players.  That is, we generalize the
standard optimization algorithm of the calculus of variations in
which functionals $F[f(x),f'(x)]$ of known form are optimized by
an independent variation of the function $f$ and the gradient
$f'$. This independent variation of each of the coordinates
$(f,f')$ over every possible value allows for instance, derivation
of the Euler-Lagrange equations providing the first order
optimization conditions. This is depicted in Fig.
\ref{f_variational_calculus} showing that every possible gradient
and trajectory through any point ``A" in the parameter space must
be considered to locate optimal trajectories. Any restrictions on
this search of all possible trajectories constrains the
optimization.  For instance, when players are restricted to using
only a particular type of probability measure space, i.e. mixed or
behavioural strategy spaces, then expected payoff functions have
fixed functional form, are continuous, and possess a single
gradient at every point in the joint function space.  These
restrictions allow use of the calculus (effectively fixed point
theorems) rather than a generalized calculus of variations to
locate equilibria.

We argue that, under CKR, players should potentially benefit from
the ability to search an enlarged mathematical space including
many alternative joint probability measure spaces. A complete
search of this enlarged mathematical space requires that they
examine not only every possible value of the expected payoff
functions at every point in their parameter space, but also every
possible gradient at every one of those points. In the following,
we show that different probability measure spaces can associate
different gradients with the same point in the joint expected
payoff function space, and we argue that every such possible
gradient must be taken into account in any complete variational
and functional optimization. That is, when players $X$ and $Y$ are
seeking to optimize their respective expected payoffs
$\langle\Pi^X\rangle$ and $\langle\Pi^Y\rangle$, they must examine
not only every possible pair of joint values
$\left(\langle\Pi^X\rangle,\langle\Pi^Y\rangle\right)$ but also
every possible joint gradient
$\left(\frac{\partial\langle\Pi^X\rangle}{\partial
p_1},\frac{\partial\langle\Pi^X\rangle}{\partial
p_2},\dots,\frac{\partial\langle\Pi^Y\rangle}{\partial
q_1},\frac{\partial\langle\Pi^Y\rangle}{\partial q_2},\dots
\right)$ evaluated with respect to every possible parameterization
$(p_1,p_2,\dots,q_1,q_2,\dots)$ defined in every possible joint
probability measure space.

\begin{figure}[htb]
\centering
\includegraphics[width=\columnwidth,clip]{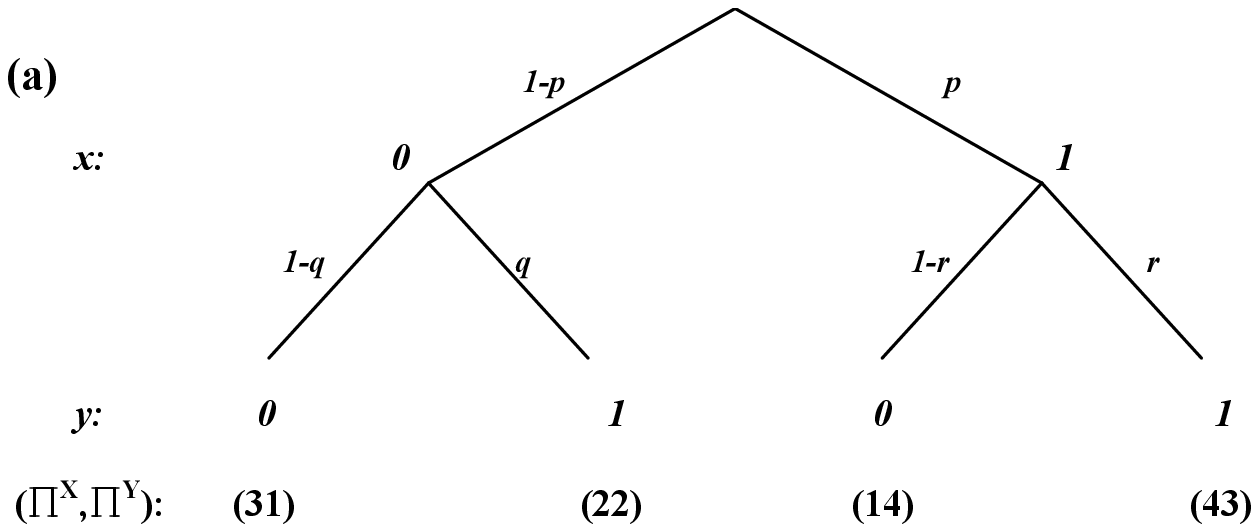}
\includegraphics[width=0.8\columnwidth,clip]{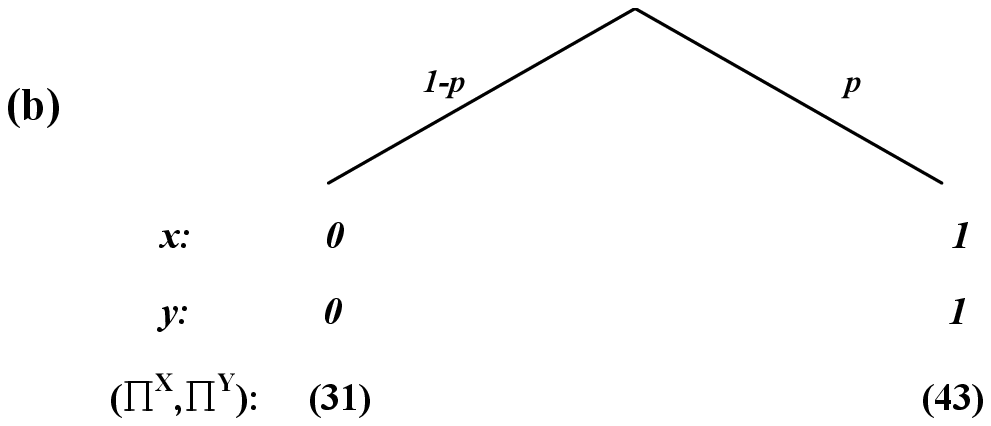}
\caption{\em  A simple example two-player-two-stage perfect
information game showing the different game decision trees
specific to players adopting the respective joint probability
spaces (a) ${\cal P}^X_0\times{\cal P}^Y_0$ where optimal player
choices are $(x,y)=(0,1)$ generating expected payoffs
$\left(\langle\Pi^X\rangle,\langle\Pi^Y\rangle\right)=(2,2)$, and
(b) ${\cal P}^X_0\times{\cal P}^Y_1$ where optimal player choices
are $(x,y)=(1,1)$ generating expected payoffs
$\left(\langle\Pi^X\rangle,\langle\Pi^Y\rangle\right)=(4,3)$.
 \label{f_simple_conventional_correlated}}
\end{figure}

\begin{figure*}[htb]
\centering
\includegraphics[width=\textwidth,clip]{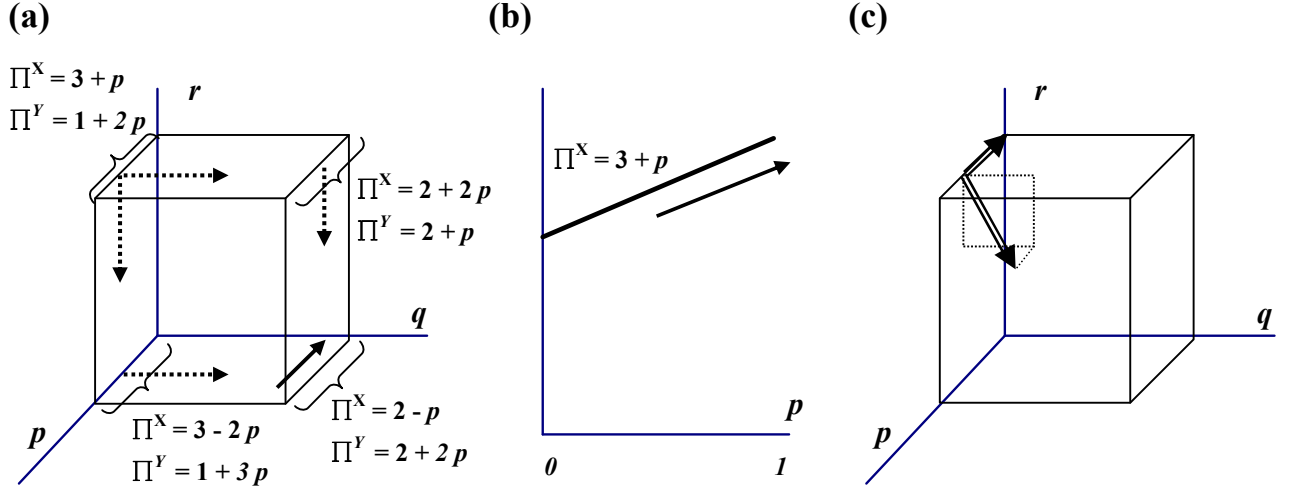}
\caption{\em (a) Game theory adopts a joint probability measure
space ${\cal P}^X_0\times{\cal P}^Y_0$ in which expected payoffs
vary over three dimensions $(p,q,r)$ and where positive gradients
with respect to $q$ and $r$ (dotted arrows) and with respect to
$p$ (solid arrow) ensure that players maximize joint payoffs by
choosing $(p,q,r)=(0,1,0)$.  (b) An alternate joint probability
space ${\cal P}^X_0\times{\cal P}^Y_1$ in which expected payoffs
vary solely over a single dimension $p$ with positive gradients
with respect to $p$ (solid arrow) ensuring that players optimize
payoffs by choosing $p=1$. (c) The choice of two alternate
probability spaces (more are possible) associates two different
total gradients (double-lined arrows) with any point along the
perfect correlation line $\rho_{xy}=1$ at $(q,r)=(0,1)$ with
$\left.\left(\frac{\partial\langle\Pi^X_{00}(p,q,r)\rangle}{\partial
p}, \frac{\partial\langle\Pi^Y_{00}(p,q,r)\rangle}{\partial q},
\frac{\partial\langle\Pi^Y_{00}(p,q,r)\rangle}{\partial
r}\right)\right|_{(q,r)\rightarrow(0,1)}\neq\left(\frac{\partial\langle
\Pi^X_{01}(p)\rangle}{\partial p}\right)$. In the absence of any
effective decision procedure privileging any one space over
another, players should examine all possible spaces, all possible
gradients, and all possible trajectories.
 \label{f_multiple_spaces}}
\end{figure*}

\section{Variational optimization in multiple stage games}

In this section, we use a simple two-player-two-stage game to
introduce standard mathematical methods that have not previously
been applied in game analysis.  Our goal is to demonstrate that
even simple games can exhibit expected payoffs with multiple
functional forms, multiple gradients and multiple trajectories at
the same point in the parameter space necessitating use of
variational optimization methods.

Suppose two players, denoted $Z\in\{X,Y\}$, each seek to optimize
their own outcomes in a strategic interaction, where in stage one
$X$ makes a choice of $x=0$ or $x=1$.  In the second stage, player
$Y$ is aware of the opponent's previous choice and must also make
a choice of $y=0$ or $y=1$ at which point the game terminates and
players obtain payoffs $\Pi^Z(x,y)$ as shown in Fig.
\ref{f_simple_conventional_correlated}(a).

Game theory analyzes this game by having each player adopt a joint
probability space allowing the complete analysis of every possible
choice that might be made in the game. In this game, both players
suppose that $X$ adopts a probability space ${\cal P}^X_0$ with
random variable $x\in\{0,1\}$ taking value $x=1$ with probability
$p$. In turn, both players suppose that player $Y$ chooses a
probability space allowing for any degree of correlation $\rho$
between the observable game events $x$ and $y$, that is, that
these variables might be perfectly correlated $\rho_{xy}=1$, or
perfectly anti-correlated $\rho_{xy}=-1$, or entirely uncorrelated
$\rho_{xy}=0$, or any value in between.  Player $Y$ does this by
adopting the probability space ${\cal P}^Y_0$ with random
variables $y,u,v\in\{0,1\}$ with $u$ and $v$ independent and
taking values $u=1$ with probability $q$ and $v=1$ with
probability $r$.  The random variable $y$ is functional determined
to be
\begin{equation}
    y = \left\{
\begin{array}{cc}
  u & \mbox{if } x=0 \\
    &   \\
  v &  \mbox{if } x=1, \\
\end{array}
    \right.
\end{equation}
giving
\begin{equation}       \label{eq_y_conditioned_on_x}
    P^Y(y|x) = \left\{
 \begin{array}{l}
  P^Y(0|0) = 1-q \\
   \\
  P^Y(1|0) = q \\
   \\
  P^Y(0|1) = 1-r \\
   \\
  P^Y(1|1) = r \\
 \end{array}
    \right. .
\end{equation}
As desired, this choice of probability space allows the players to
examine every possible correlation state between $x$ and $y$
defined as
\begin{eqnarray}
    \rho_{xy}(p,q,r)
    &=& \frac{\langle xy\rangle-\langle x\rangle\langle y\rangle}{
    \sqrt{\langle x^2\rangle-\langle x\rangle^2}
    \sqrt{\langle y^2\rangle-\langle y\rangle^2}} \nonumber \\
    &=& \frac{\sqrt{p(1-p)}(r-q)}{
    \sqrt{\left[q+p(r-q)\right]\left[1-q-p(r-q)\right]}}.
\end{eqnarray}
Then, $x$ and $y$ are perfectly correlated at
$\rho_{xy}(p,0,1)=1$, perfectly anti-correlated at
$\rho_{xy}(p,1,0)=-1$, and uncorrelated if either $p=0$ or $1$ or
$q=r$ giving $\rho_{xy}=0$.  As shown in Fig.
\ref{f_multiple_spaces}(a), in the joint probability space ${\cal
P}^X_0\times{\cal P}^Y_0$, the expected payoff functions are
\begin{eqnarray}
 \langle\Pi^X_{00}(p,q,r)\rangle
      &=& \sum_{xy=0}^1 P^X(x) P^Y(y|x) \Pi^X(x,y)  \nonumber \\
      &=& 3-q+p(q+3r-2)  \nonumber \\
 \langle\Pi^Y_{00}(p,q,r)\rangle
      &=& \sum_{xy=0}^1 P^X(x) P^Y(y|x) \Pi^Y(x,y)  \nonumber \\
      &=& 1+q-p(q+r-3),
\end{eqnarray}
so the gradients with respect to the three continuous dependent
variables $p$, $q$ and $r$ are
\begin{eqnarray}
    \frac{\partial\langle \Pi^X_{00}(p,q,r)\rangle}{\partial p}
      &=& q + 3 r -2 \nonumber \\
    \frac{\partial\langle \Pi^Y_{00}(p,q,r)\rangle}{\partial q}
      &=&  1-p \nonumber \\
    \frac{\partial\langle \Pi^Y_{00}(p,q,r)\rangle}{\partial r}
      &=&  -p.
\end{eqnarray}
As shown in Fig. \ref{f_multiple_spaces}(a), this
three-dimensional gradient exists and is non-zero even when $x$
and $y$ are perfectly correlated $\rho_{xy}=1$ at all points
$(q,r)=(0,1)$ so payoffs are not optimized at these points.  In
fact, given the choice of probability space ${\cal
P}^X_0\times{\cal P}^Y_0$, both players conclude that $Y$
maximizes their payoff by setting $(q,r)=(1,0)$ while $X$
maximizes their payoff by setting $p=0$. The resulting move
choices are $(x,y)=(0,1)$ generating payoffs of
$\left(\langle\Pi^X\rangle,\langle\Pi^Y\rangle\right)=(2,2)$. This
completes our analysis of the usually adopted joint probability
measure space, and we now turn to examine alternatives.

In any game, alternate joint probability measure spaces exist with
expected payoff functions of different functional form and
different gradients at the same point in the parameter space.
Suppose that player $Y$ chooses a different probability space
${\cal P}^Y_1$ in which they treat the observed value of the
random variable $x$ as a coin toss determining their choice of
$y=1$ with probability $p$. That is, $Y$ functionally assigns the
random variable $y$ to be perfectly correlated with the observed
random variable $x$ via
\begin{eqnarray}
   y &=& x \nonumber \\
   P^Y(y|x) &=& \delta_{yx}.
\end{eqnarray}
This functional assignment does not require any communication
between player $X$ and $Y$. Then, in the joint probability space
${\cal P}^X_0\times{\cal P}^Y_1$, the expected payoff functions
straightforwardly equal
\begin{eqnarray}
 \langle\Pi^X_{01}(p)\rangle
      &=& \sum_{xy=0}^1 P^X(x) P^Y(y|x) \Pi^X(x,y)  \nonumber \\
      &=&   3 + p \nonumber \\
 \langle\Pi^Y_{01}(p)\rangle
      &=& \sum_{xy=0}^1 P^X(x) P^Y(y|x) \Pi^Y(x,y)  \nonumber \\
      &=& 1+ 2 p,
\end{eqnarray}
as seen in the decision tree of Fig.
\ref{f_simple_conventional_correlated}(b), and in the expected
payoff function space of Fig. \ref{f_multiple_spaces}(b). These
expected payoff functions are now dependent only on the single
freely varying parameter $p$ determining the gradient with respect
to $p$ to be
\begin{equation}
    \frac{\partial\langle \Pi^X_{01}(p)\rangle}{\partial p} = 1.
\end{equation}
Consequently, player $X$ maximizes their payoff by setting $p=1$
to choose $x=1$ leading $Y$ to set $y=1$.  Thus, in the joint
probability space ${\cal P}^X_0\times{\cal P}^Y_1$, player payoffs
are $\left(\langle\Pi^X\rangle,\langle\Pi^Y\rangle\right)=(4,3)$.

\begin{figure*}[htb]
\centering
\includegraphics[width=0.8\textwidth,clip]{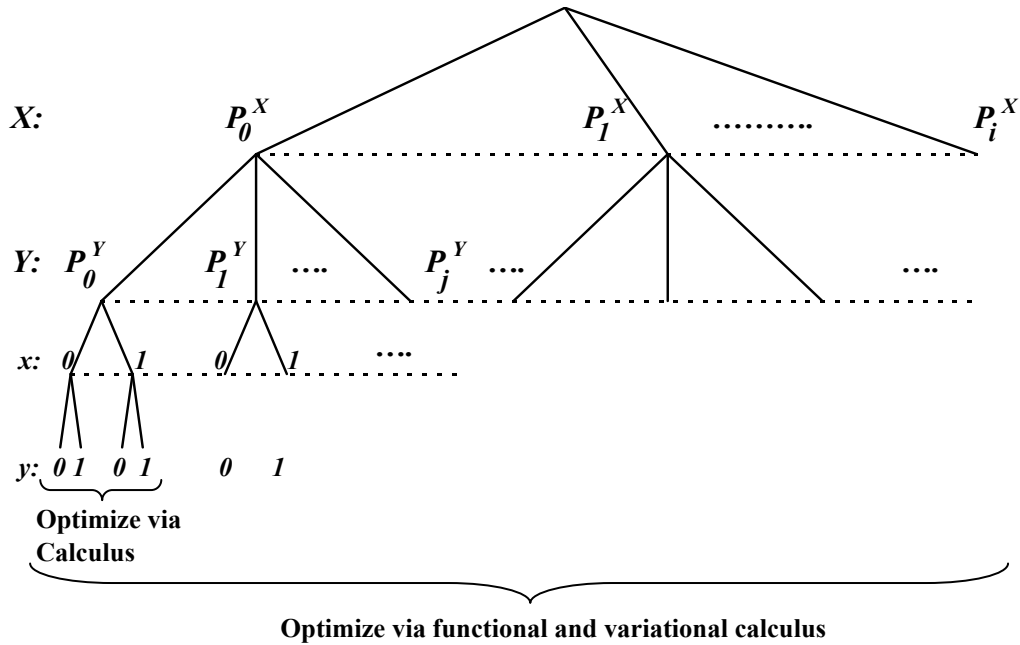}
\caption{\em  A schematic representation of a complete analysis of
the simple game of Fig.
\protect\ref{f_simple_conventional_correlated} showing that each
player must first decide which probability space to adopt. Here,
players $X$ and $Y$ simultaneously choose their respective
probability spaces ${\cal P}^X_i$ and ${\cal P}^Y_j$ from among an
infinite number of possibilities, where each choice generates a
different decision tree for the game defining altered payoff
functions.  Players do not know their opponent's choice of
probability space. Given the adoption of a particular joint
probability space ${\cal P}^X_i\times{\cal P}^Y_j$, expected
payoff functions are known and continuous in terms of their input
parameters so the calculus, suitably generalized for
multiple-player strategic interactions, can be used to optimize
payoffs. However, optimizing the choice of joint probability space
involves expected payoff functionals requiring players to use
variational calculus methods to optimize payoffs.
 \label{f_total_tree}}
\end{figure*}

We now have two possible joint probability spaces; that normally
adopted in game theory ${\cal P}^X_0\times{\cal P}^Y_0$ and the
novel ${\cal P}^X_0\times{\cal P}^Y_1$.  In these alternate
spaces, the expected payoff functions possess exactly the same
value when $x$ and $y$ are perfectly correlated but possess
entirely different gradients at this point---see Fig.
\ref{f_multiple_spaces}(c). Variational optimization principles
insist that every possible functional form and gradient must be
taken into account in any complete optimization.  These principles
permit players to infinitely vary the ``immutable" functional
assignments defining any space (i.e. $y=\delta_{x0}u+\delta_{x1}v$
and $y=x$ above), providing access to a vastly larger decision
space than usually analyzed in game theory. It is not a question
of which space is best, rather, it is a question of either
restricting the analysis to a single space or allowing players to
analyze all possible spaces.

Game theory adopts expected payoff ``functions"
$\langle\Pi^Z_{00}(p,q,r)\rangle$ allowing examination of every
possible combination of payoff values and assumes that this is
sufficient for optimization.  However, while these functions can
duplicate every possible payoff value, they cannot duplicate every
possible functional gradient---and optimization depends on
gradients. When $Y$ adopts a randomization device (a ``roulette")
which perfectly correlates $x$ and $y$ via the probability space
${\cal P}^Y_1$, then certainly
$\langle\Pi^Z_{00}(p,0,1)\rangle=\langle\Pi^Z_{01}(p)\rangle$, but
these functions have different dimensionality and gradients.  That
is,
$\left.\left(\frac{\partial\langle\Pi^X_{00}(p,q,r)\rangle}{\partial
p}, \frac{\partial\langle\Pi^Y_{00}(p,q,r)\rangle}{\partial q},
\frac{\partial\langle\Pi^Y_{00}(p,q,r)\rangle}{\partial
r}\right)\right|_{(q,r)\rightarrow(0,1)}\neq\left(\frac{\partial\langle
\Pi^X_{01}(p)\rangle}{\partial p}\right)$.  Similar results apply
for points at different correlation values $\rho_{xy}$; should $Y$
adopt a randomization device where $y$ is entirely uncorrelated
with $x$ via a new probability space ${\cal P}^Y_2$, then
certainly
$\langle\Pi^Z_{00}(p,q,q)\rangle=\langle\Pi^Z_{02}(p,q)\rangle$
but these functions again have different dimensionality and
gradients
$\left.\left(\frac{\partial\langle\Pi^X_{00}(p,q,r)\rangle}{\partial
p}, \frac{\partial\langle\Pi^Y_{00}(p,q,r)\rangle}{\partial q},
\frac{\partial\langle\Pi^Y_{00}(p,q,r)\rangle}{\partial
r}\right)\right|_{r=q}\neq\left(\frac{\partial\langle
\Pi^X_{02}(p,q)\rangle}{\partial p},\frac{\partial\langle
\Pi^Y_{02}(p,q)\rangle}{\partial q}\right)$.  These inequalities
result as the usually adopted space ${\cal P}^X_0\times{\cal
P}^Y_0$ evaluates gradients using infinitesimals between points
with different correlations so
$\Delta\rho_{xy}=\rho_{xy}(p,q,r)-\rho_{xy}(p+\delta p,q+\delta
q,r+\delta r)\neq 0$.  In contrast, when a roulette possesses a
known correlation state as in the spaces ${\cal P}^X_0\times{\cal
P}^Y_1$ or ${\cal P}^X_0\times{\cal P}^Y_2$, then gradients are
evaluated taking all constraints into account ensuring
$\Delta\rho_{xy}=0$. Game analysis does not include every possible
correlation constraint or every possible roulette, and taking
these alternatives into account requires the variational methods
presented in this paper.

\begin{figure*}[htb]
\centering
\includegraphics[width=0.8\textwidth,clip]{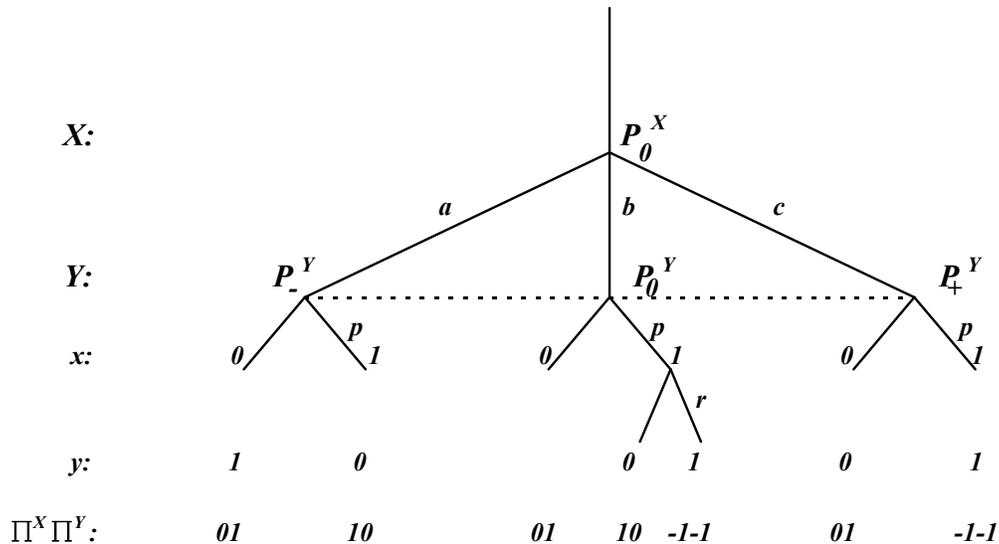}
\caption{\em An extended decision tree for the chain store paradox
where player $X$ adopts the usual probability space ${\cal P}^X_0$
with certainty while $Y$ has a choice of three alternate
probability spaces ${\cal P}^Y_j$ for $j\in\{-,0,+\}$ selected
with probabilities $a$, $b$, and $c$, and respectively denoting
anti-correlated, independent, and correlated decision making. Game
theory examines only the branch generated by the joint probability
space ${\cal P}^X_0\times{\cal P}^Y_0$ in which a potential new
market entrant $X$ must decide to either stay out of a new market
$x=0$ or enter the market $x=1$ (with probability $p$), in which
case the monopolist $Y$ chooses to either acquiesce $y=0$ or fight
their entry $y=1$ (with probability $r$), with the corresponding
payoffs shown. \label{f_chain_total_tree}}
\end{figure*}

We suggest that the example game described in the two decision
trees of Fig. \ref{f_simple_conventional_correlated} is best seen
as having the schematic form shown in Fig. \ref{f_total_tree} in
which players must first choose which probability space they will
adopt, where this choice is unknown to their opponents at the
commencement of the game, and must then optimize their payoffs
given the possible joint probability spaces that might be adopted.
In such generalized trees incorporating choice of probability
space, standard approaches can be applied to locate pure
``variational" strategies, probabilistic ``variational" mixed and
``variational" behavioural strategies, and ``variational"
equilibria.  Of course, introducing ``variational" mixed and
behavioural strategies means that players must introduce yet
further probability spaces allowing the optimization of these
probabilistic strategies.

To provide a concrete illustration of our approach, we now show
that rational players using variational optimization methods can
resolve the chain store paradox.

\section{Resolving the chain store paradox}

A minimal chain store paradox is shown as the central branch in
Fig. \ref{f_chain_total_tree} generated by the adoption of the
joint probability space ${\cal P}^X_0\times{\cal P}^Y_0$ (defined
below).  This game is played over two sequential stages where
first, a potential market entrant $X$ must decide to either stay
out of a new market $x=0$ or enter that market $x=1$. Their
opponent, the monopolist $Y$, observes this choice. Should no
market entry occur, $X$ neither gains nor loses any payoff while
$Y$ gains monopolist profits so $(\Pi^X,\Pi^Y)=(0,1)$. In
contrast, should $X$ enter the market, $Y$ must then decide
whether to acquiesce to their opponent's entry $y=0$ by leaving
prices unchanged and sharing profits so $(\Pi^X,\Pi^Y)=(1,0)$, or
by driving $X$ out of business by price cutting so payoffs are
$(\Pi^X,\Pi^Y)=(-1,-1)$.

A backwards induction analysis of the central branch of Fig.
\ref{f_chain_total_tree} in isolation indicates that $X$ will
enter the market confident that the monopolist will not forego
profits to fight their entry \cite{Selten_78_12}.  Based on this,
many economists argue it is irrational for monopolists to engage
in predatory pricing to drive rivals out of business as predation
is costly while potential new entrants well understand that price
cuts are temporary and monopoly profits readily attract new market
entrants \cite{Milgrom_1982_280}. Efforts to resolve the paradox
include introducing multiple stages permitting reputation and
deterrence effects \cite{Selten_78_12}, as well as asymmetric
information, mistakes, bounded rationality or imperfect
information and uncertainty
\cite{Milgrom_1982_280,Rosenthal_1981_92,Davis_85_13,Kreps_1982_253,Trockel_86_16}.
For a review, see \cite{Wilson_1992_305}.

The decision and payoff combinations above define the general
strategic optimization problem faced by the players in the chain
store paradox as
\begin{eqnarray}
    X\!\!: \max_{{\cal P}^X}\;\; \langle\Pi^X\rangle
      &=& P^{X}(1) \left[1 - 2 P^Y(1|1) \right] \nonumber  \\
      && \\
    Y\!\!: \max_{{\cal P}^Y}\;\; \langle\Pi^Y\rangle
      &=& 1-P^{X}(1)  - P^{X}(1)P^Y(1|1). \nonumber
\end{eqnarray}
Here, players alter their choice of joint probability space ${\cal
P}^X\times{\cal P}^Y$ such that, within the selected optimal joint
probability space, the optimization of their respective
probability distributions $P^X(x)$ and $P^Y(y|x)$ allows optimal
choices $(x,y)$ to be made so as to maximize respective payoffs.

A complete derivation of the ``variational" equilibria of the
extended tree of Fig. \ref{f_chain_total_tree} is of course
possible (and indicates that as long as $Y$ sets
$c+br\geq\frac{1}{2}$, most easily achieved by choosing $c=1$,
then $X$ maximizes their payoff through the choice $p=0$). In this
paper, we locate pure variational equilibria in the chain store
game. That is, we suppose that player $X$ always chooses ${\cal
P}^X_0$ (other choices are possible) while player $Y$ chooses with
certainty any of the three probability spaces ${\cal P}^Y_j$ with
$j\in\{-,0,+\}$.  The interpretation of these latter spaces is
that $j=$``$-$" indicates that $y$ is perfectly anti-correlated to
$x$, $j=$``$0$" indicates that $y$ is entirely independent, and
$j=$``$+$" indicates that $y$ is perfectly correlated to $x$---see
Fig. \ref{f_chain_total_tree}.

First, we replicate the usual game analysis by supposing that
player $X$ adopts a probability space ${\cal P}^X_0$ with random
variable $x\in\{0,1\}$ such that $x=1$ with probability $p$, while
player $Y$ adopts the probability space ${\cal P}^Y_0$ with random
variables $y,u,v\in\{0,1\}$ such that $u$ and $v$ are independent
random variables taking value $u=1$ with probability $q$ and $v=1$
with probability $r$, and where the random variable $y$ is
functionally assigned as $y=\delta_{x0}u+\delta_{x1}v$.
Altogether, this gives the same probability parameterization for
$P^Y(y|x)$ as appears in Eq. \ref{eq_y_conditioned_on_x}. In the
joint probability space ${\cal P}^X_0\times{\cal P}^Y_0$, the
optimization problem reduces to
\begin{eqnarray}
    X\!\!: \max_{p}\;\; \langle\Pi^X_{00}\rangle
      &=& P^{X}(1) \left[1 - 2 P^Y(1|1) \right] \nonumber \\
      &=& p ( 1 - 2 r ) \nonumber \\
      && \\
    Y\!\!: \max_{r}\;\; \langle\Pi^Y_{00}\rangle
      &=& 1-P^{X}(1)  - P^{X}(1)P^Y(1|1) \nonumber \\
      &=& 1-p  - pr,  \nonumber
\end{eqnarray}
so the gradients with respect to the two continuous dependent
variables $p$ and $r$ are
\begin{eqnarray}
    \frac{\partial\langle\Pi^X_{00}\rangle}{\partial p}
    &=& 1 - 2 r  \nonumber  \\
    \frac{\partial\langle\Pi^Y_{00}\rangle}{\partial r}
    &=& - p.
\end{eqnarray}
Essentially then, the monopolist $Y$ maximizes their expected
payoff by setting $r=0$ and always acquiesces to new market
entrants, while $X$ maximizes their payoff by choosing $p=1$ and
so always decides to enter the market. The resulting expected
payoffs given that players adopt this sole perfect Nash equilibria
of $(p,r)=(1,0)$ are
$\left(\langle\Pi^X\rangle,\langle\Pi^Y\rangle \right)=(1,0)$.

Suppose however that players $X$ and $Y$ choose the joint
probability space ${\cal P}^X_0\times{\cal P}^Y_+$ (the rightmost
branch of Fig. \ref{f_chain_total_tree}) where $y$ is perfectly
correlated with $x$ via the functional assignment $y=x$ and
$P(y|x)=\delta_{yx}$ altering the strategic optimization problem
to
\begin{eqnarray}
    X\!\!: \max_{p}\;\; \langle\Pi^X_{0+}\rangle
      &=& P^{X}(1) \left[1 - 2 P^Y(1|1) \right] \nonumber \\
      &=& - p \nonumber \\
      &&  \\
    Y\!\!:  \langle\Pi^Y_{0+}\rangle
      &=& 1-P^{X}(1)  - P^{X}(1)P^Y(1|1) \nonumber \\
      &=& 1-2p.  \nonumber
\end{eqnarray}
These expected payoff functions are continuous over the single
freely varying parameter $p$ giving the gradient
\begin{equation}
    \frac{\partial\langle\Pi^X_{0+}\rangle}{\partial p} =
      -1 < 0
\end{equation}
ensuring that player $X$ maximizes their expected payoff by
setting $p=0$ and not entering the market. That is, when players
$(X,Y)$ adopt the ${\cal P}^X_0\times{\cal P}^Y_+$ joint
probability space, they maximize their payoffs via the combination
$(x,y)=(0,0)$ to garner payoffs
$\left(\langle\Pi^X\rangle,\langle\Pi^Y\rangle\right)= (0,1)$.

Alternatively, in the anti-correlated joint probability space
${\cal P}^X_0\times{\cal P}^Y_-$ (the leftmost branch of Fig.
\ref{f_chain_total_tree}), $y$ is perfectly anti-correlated with
$x$ via $y=(1-x)$ and $P(y|x)=\delta_{y(1-x)}$ giving the altered
strategic optimization problem of
\begin{eqnarray}
    X\!\!: \max_{p}\;\; \langle\Pi^X_{0-}\rangle
      &=& P^{X}(1) \left[1 - 2 P^Y(1|1) \right] \nonumber \\
      &=& p \nonumber \\
      &&  \\
    Y\!\!:  \langle\Pi^Y_{0-}\rangle
      &=& 1-P^{X}(1)  - P^{X}(1)P^Y(1|1) \nonumber \\
      &=& 1-p.  \nonumber
\end{eqnarray}
Again, these are functions of the sole parameter $p$ giving the
gradients \begin{equation}
    \frac{\partial\langle\Pi^X_{0-}\rangle}{\partial p} =
      1 > 0,
\end{equation}
ensuring that player $X$ sets $p=1$ and chooses to enter the
market.  The result is that when players $(X,Y)$ adopt the ${\cal
P}^X_0\times{\cal P}^Y_-$ joint probability space, they maximize
their payoffs via the combination $(x,y)=(1,0)$ to garner payoffs
$\left(\langle\Pi^X\rangle,\langle\Pi^Y\rangle\right)= (1,0)$.

Altogether, when players consider only pure variational strategies
(specifying probability spaces and move choices), the various
payoffs available are
\begin{equation}
    \begin{array}{c|c}
 \left(\langle\Pi^X\rangle,\langle\Pi^Y\rangle\right) & {\cal P}^X_0 \\ \hline
                                                &               \\
       {\cal P}^Y_-                             &      (1,0)    \\
                                                &               \\
       {\cal P}^Y_0                             &      (1,0)    \\
                                                &               \\
       {\cal P}^Y_+                             &      (0,1),    \\
    \end{array}
\end{equation}
making it evident that to maximize their expected payoff, player
$Y$ must rationally elect to use probability space ${\cal P}^Y_+$
in preference to either ${\cal P}^Y_0$ or ${\cal P}^Y_-$. That is,
$Y$ will undertake to functionally correlate their move to the
previous choice of the potential market entrant, and thereby deny
themselves a choice about the setting of $y$ once the game has
commenced.  In the probability space ${\cal P}^Y_+$, the
optimization by player $Y$ has no second stage component as the
joint probability distributions are inseparable, and an
opportunity for a second stage optimization exist only in the
space ${\cal P}^Y_0$. Player $Y$ foregoes a choice during the game
itself knowing this to be payoff maximizing.  Player $X$, being
aware of this will not enter the market as in the minimal chain
store game described by ${\cal P}^X_0\times{\cal P}^Y_+$, entry
automatically invokes retaliation. We thus reconcile game
theoretic prediction and observed human behaviour implying human
players generally commence a strategic analysis by first
optimizing their choice of probability space and only subsequently
optimizing the probability distributions defined by that space.

It is of course possible to consider a broader range of joint
probability spaces for both players $X$ and $Y$, though this will
not substantially alter the conclusion here that it can be
rational for a monopolist to punish market entrants to resolve the
chain store paradox.

\section{Conclusion}

This paper locates strategic equilibria using generalized calculus
of variation techniques and is thus consistent with, and extends,
the more usual methods of game theory based on the fixed point
theorems of the calculus \cite{Hart_92_19,Sorin_92_c4}. We hold
that under CKR, players might often improve their outcomes by
expanding their mathematical search space to include alternate
probability spaces. Consequently, we allow players to first
optimize their choices of probability measure space which alters
both the expected payoff functionals (not functions) and the joint
probability distributions specifying move choices to locate
``variational" equilibria. Generally, these variational equilibria
differ from Nash equilibria even in perfect information games such
as the chain store paradox considered here.  This is because
first, players are uncertain about which joint probability measure
space is in play, and second, each alternative space introduces
different correlations rendering the joint probability
distribution inseparable and altering allowable subgame
decompositions and the backwards induction analysis. We show that
when rational players variationally optimize their choice of
probability measure space to access ``variational" equilibria,
then this can reconcile game theoretic prediction and observed
human behaviour.  To illustrate this, we demonstrated that our
general variational and functional optimization approaches resolve
the chain store paradox. This strongly suggests that rational
players should, in fact, exploit unrestricted optimization in
general.

More generally, we suggest that selfish ``homo economicus" might
exploit variational optimization to access alternate probability
measure spaces to exhibit altruistic or cooperative behaviour
whenever that is payoff-maximizing.  This might explain for
instance, the efficacy of state led development processes
\cite{World_Bank_93,Stiglitz_00} and the industry wide
correlations of the ``Just-In-Time" Toyota production system
\cite{Lu_89,Womack_1990}.  Consequently, variational optimization
may change the orientation and methods of evolutionary game theory
\cite{Maynard_Smith_1982} and quantum game theory
\cite{Meyer_99_10}, among other fields, while the need to search
infinite numbers of joint probability spaces will reinforce the
importance of learning, a principal feature of evolutionary
economics \cite{Wit_1993}.  Similarly, variational optimization
will impact ``selfish gene" theory which presently holds that
genes optimize their fitness independently so altruism is
explained by relatedness and the likelihood of shared genes
\cite{Dawkins_1976}. In contrast, we suggest that modelling the
evolution of emergent hierarchical complexity in living organisms
requires taking account of alternative probability spaces
correlating system components; correlated entities together
constitute an indivisible unit which must be optimized as a whole.
Consequently, complex multicellular eukaryotes might well have
optimized fitness by adopting correlating signals to multitask
their dynamics by exploring alternate dynamical decision trees
(organismal probability spaces) most likely by expansion of their
RNA signaling capabilities
\cite{Mattick_94_823,Mattick_02_1611,Mattick_04_316}. Similar
considerations mean that neural networks can endogenously modify
correlations among their components to explore an infinite number
of alternate dynamical trees to implement complex cognition. Game
analysis also underlies the tree search ``minimax" algorithms of
artificial intelligence \cite{Hart_92_c2,Russell_2003,Callan_2003}
which typically fail to emulate human intelligence.  In chess
playing, for instance, expert human players typically employ
pattern recognition and ``chunking" \cite{Luger_1998}, and appear
to be exploiting the same correlation information that underpins
variational optimization. This is consistent with the ``social
intelligence" explanation for the runaway evolution of primate
intelligence where individuals dynamically realign their strategic
partnerships to correlate behaviours to optimize outcomes in
competitive group settings \cite{Byrne_1988,Rifkin_1995}.

It has long been thought that any strategic optimization problem
was essentially equivalent to a possibly greatly enlarged
non-strategic optimization problem. This equivalence arises as
each player can introduce sufficient new variables to fully model
all of the possible actions of all of their opponents. As a
result, strategic optimization has been thought to be of
equivalent complexity to, for instance, non-strategic physics
optimization problems, and solved by similar methods such as the
calculus of variations or the calculus.  Some have argued the
converse.  A perceived fundamental incompatibility between physics
and biological complexity motivated Mayr to claim that biology is
an autonomous science rather than a subbranch of the physical
sciences \cite{Mayr_2004}, with the factor missing in physics but
present in biology being identified as ``entailment" (essentially
correlation) by Rosen \cite{Rosen_1991}, while it has been
unconventionally argued that information science is incomplete and
that it is our growing understanding of genomic programming and
biological complexity that will contribute significant new
insights in this field (J. S. Mattick, personal communication).
Variational optimization might well help close these perceived
gaps.


\end{document}